\documentclass[12pt]{amsart}
\usepackage{amsthm,amssymb}
\newcommand{\ol}{\overline}

\newcommand{\iy}{\infty}

\newcommand{\vf}{\varphi}

\newcommand{\Adj}{\mathop{\rm Adj}}
\newcommand{\diag}{\mathop{\rm diag}}

\newcommand{\bT}{\mathbb{T}}

\newtheorem{lemma}{\bf \sc Lemma}

\theoremstyle{remark}

\renewenvironment{proof}{ P\,r\,o\,o\,f.}{$\Box$}

\begin{document}

\title[Matrix Spectral Factorization]
    {An Analytic Proof of the Matrix Spectral Factorization Theorem}

\author{L. Ephremidze, \fbox{G. Janashia}\,, and E. Lagvilava  }

\footnotetext { This paper will appear in Georgian Mathematical
Journal.}

\begin{abstract}
    An analytic proof is proposed of Wiener's theorem on
    factorization of positive definite matrix-functions.
\end{abstract}

\maketitle

\subsection{Introduction}
The series of papers of Wiener [15], [16], [18], and Helson and
Lowdenslager [4], [5]  led to the following

\smallskip

{\bf Matrix Spectral Factorization Theorem:} {\em Let
\begin{equation}
    S(z)=\big(S_{jk}(z)\big)_{j,k=1}^r\,,
\end{equation}
$|z|=1$, be a positive definite $r\!\times\!r$ matrix-function with
integrable entries, $S_{jk}(z)\in L_1(\bT)$, defined on the unit
circle. If
\begin{equation}
\log \det S(z)\in L_1({\mathbb T}),
\end{equation}
then $S(z)$ admits a factorization
\begin{equation}
    S(z)=\chi^+(z)(\chi^+(z))^*,
\end{equation}
where $\chi^+(z)$, $|z|<1$, is an analytic $r\!\times\!r$
matrix-function with entries from the Hardy space $H_2$, and
$\det\chi^+(z)$  is an outer function.}
\smallskip

The equation (3) is assumed to hold a.e. on the unit circle $\bT$
and $(\chi^+)^*=(\ol{\chi^+})^T$ is the Hermitian conjugate of
$\chi^+$.

The spectral factorization (3) is unique up to a constant unitary
multiplier, and the unique spectral factor with additional
requirement that $\chi^+(0)$ is positive definite is called
canonical.

It is well-known that the condition (2) is also necessary for the
existence of the factorization (3).

Spectral factorization problem was originally formulated in the
scalar case by Wiener [14] and Kolmogorov [9], independently from
each other, in connection with developed linear prediction theory of
stationary stochastic processes. Wiener [16] used the matrix
generalization of the problem in the study of multidimensional
stationary time series. Since then the spectral factorization has
became an important tool in solution of various applied problems in
Control Engineering and Communications.

The prediction theory of stationary processes itself, namely the
idea derived from the theory of linear least-squares method, is
helpful to prove the spectral factorization theorem in the scalar as
well in the matrix case and Wiener successfully exploited this idea.
Helson and Lowdenslager derived the theorem from the theory of
invariant subspaces and this proof, together with its
generalizations to positive  operator valued functions, is well
presented in [3].

However, in the scalar case, the theory of the Hardy spaces $H_p$
provides a simple solution to the spectral factorization problem.
Furthermore, the spectral factor can be written in an explicit form
(see [10])
\begin{equation}
\chi^+(z)=\exp\left(\frac 1{4\pi} \int_0^{2\pi}
\frac{e^{it}+z}{e^{it}-z}\log S(e^{it})\,dt\right).
\end{equation}
No analog of the formula (4) exists in the matrix case since the
equation $\exp(A+B)=\exp A\exp B$ breaks down for non-commutative
matrices and, in general, the factorization problem becomes more
demanding. Nevertheless, one can observe that the formulation of the
spectral factorization theorem in the matrix case does not go beyond
the theory of $H_p$ spaces and naturally its proof had been searched
within this theory. In the present paper we provide such proof on 50
years anniversary of this important theorem.

Together with the existence theorem of the spectral factorization,
the challenging problem became the actual approximate computation of
the spectral factor $\chi^+(z)$ for a given matrix-function (1)
since, as it was mentioned above, this procedure found its immediate
applications in solutions of various practical problems. Wiener [17]
made first efforts to create a sound computational algorithm of
matrix spectral factorization and since then several dozens of
papers have addressed to the development of such algorithm. It is
surely hopeless to give here a short but still comprehensive account
on all these developments (see [7], [13]). Wiener and his followers
[17], [11] as well as Kolmogorov's students in Russia [12], [20]
were using methods of Functional Analysis for solution of the
problem. A Newton-Raphson type iterative algorithm was constructed
by G. Wilson [19]. Kalman's state-space approach to the Wiener
filtering problem [8] became fruitful for spectral factorization as
well, and by this way the spectral factorization problem was
actually reduced to the solution of algebraic Riccati equation.
Until recently this method was considered as most effective for
matrix spectral factorization (see [7], p. 206), though it can be
used only for rational matrix functions. Numerous algorithmic type
improvements were proposed by the different authors in this
direction.

An absolutely new approach to the spectral factorization problem in
the matrix case was developed in [6] (see also [1]). Without
imposing on matrix-function (1) any additional restriction, apart
from the necessary and sufficient condition (2), effective method of
approximate computation of $\chi^+$ is proposed. This is the first
time that the methods of Complex Analysis and Hardy Spaces were used
for construction of the algorithm, which naturally turned out to be
very efficient as corresponding software implementation confirms
(see demo version at www.rmi.acnet.ge/SpFact). The decisive role of
unitary matrix-functions in the factorization process is revealed
which, by flexible manipulations, completely absorbs all the
technical difficulties of the problem leaving very few and simple
procedures for computation. This method of matrix spectral
factorization is undoubtedly as simple as possible and in the
present paper we would like to demonstrate that in fact it contains
an analytic proof of the existence theorem as well, which was
formulated in the beginning. No other assumptions beyond the theory
of Hardy spaces are used for this purpose.

\subsection{Notation}
$L_p(\bT)$, $p>0$, is the class of $p$-integrable complex functions
defined on the unit circle, and $H_p$ denotes the Hardy space of
analytic functions in the unit disk,
$$
H_p=\left\{f\in {\mathcal A}(D):\sup_{r<1}
\int\nolimits_0^{2\pi}|f(re^{it})|^p\,dt<\iy\right\}.
$$
A function $Q(z)\in H_p$ is called outer, we denote $Q(z)\in H_p^O$,
if
$$
Q(z)=c\cdot\exp\left(\frac 1{2\pi}
\int\nolimits_0^{2\pi}\frac{e^{it}+z}{e^{it}-z}\log
|Q(e^{it})|\,dt\right),
$$
where $|c|=1$.

The $n$th Fourier coefficient of $f\in L_1(\bT)$ is denoted by
$c_n(f)$ and, for $p\geq 1$, $L_p^+({\mathbb T}):=\{f\in
L_p({\mathbb T}): c_n(f)=0 \text{ for } n<0\}$. The spaces
$L_p^+({\mathbb T})$ and $H_p$ are naturally identified, so that we
can speak about the value of function $f\in L_p^+({\mathbb T})$ in
$z\in D$.

Let $L_p^{N_-}(\bT)$, $N>0$, $p\geq 1$, be the set of functions $f$
from $L_p(\bT)$ for which $c_n(f)=0$ whenever $n<-N$, and let
$L_{[N]}^+$
 be the set of analytic polynomials whose nonzero
coefficients range from $0$ to $N$.

If $M $ is a matrix (resp. matrix-function), then $\ol{M}$ denotes
the matrix (resp. matrix-function) with conjugate entries and
$M^*:=\ol{M}^T$. The upper-left $m\!\times\!m$ submatrix of $M$ is
denoted by $M^{[m,m]}$.

We say that matrix-functions have some property, say, belong to
$L_p({\mathbb T})$ or are convergent, etc, if their entries have
this property.

A $r\!\times r$ matrix $U$ is called unitary if $ UU^*=I_r $, where
$I_r$ is the identity matrix of dimension $r$. Since the rows and
columns of $U$ are orthonormal, its entries are bounded by $1$. A
unitary matrix-function $U(z)$ means that it is unitary for a.a.
$z\in\bT$.

The notation $\diag(u_1,u_2,\ldots,u_r)$ stands for the diagonal
$r\!\times\!r$ matrix with corresponding entries on the main
diagonal.

\subsection{The uniqueness of spectral factorization} In the proof
of the convergence property of the above mentioned algorithm given
in [6], at least formally, the existence of spectral factorization
is used. However, one can observe that rather the uniqueness than
the existence of spectral factorization provides this convergence.
So we start with a simple proof of the uniqueness theorem (cf.
 [2]), emphasizing that it can be obtained without a priori knowledge
of the existence of the spectral factorization itself.

We use the following generalization of Smirnov's theorem concerning
functions from the Hardy spaces $H_p$ (see [10], p.~109): {\em Let
$f(z)=g(z)/h(z)$, where $g\in H_{p_1}$, $p_1>0$, and $h\in
H^O_{p_2}$, $p_2>0$. If the boundary values $f(e^{it})\in
L_p({\mathbb T})$, $p>0$, then $f\in H_p$\,.}

\smallskip

{\bf Uniqueness Theorem:} {\em If
\begin{equation}
    S(z)=\chi_j^+(z)(\chi_j^+(z))^*,\;\;\;\;\;j=1,2,
\end{equation}
are two spectral factorizations of a given spectral density $S(z)$,
then
\begin{equation}
    \chi_1^+(z)=\chi_2^+(z)\cdot U,\;\;\;\;\;\;|z|<1,
\end{equation}
for some constant unitary matrix $U$.}
\smallskip

\begin{proof}
It follows from (5) that $\chi_1^+(z)(\chi_1^+(z))^*=
\chi_2^+(z)(\chi_2^+(z))^*$ a.e. on $\bT$, so that
\begin{equation}
 (\chi_2^+(z))^{-1}\chi^+_1(z)\big((\chi_2^+(z))^{-1}\chi^+_1(z)\big)^*=I_r
\text{ for a.a. }z\in \bT.
\end{equation}
Thus the analytic matrix-function
\begin{equation} U(z):=(\chi_2^+(z))^{-1}\chi^+_1(z),\;\;\;|z|<1,\end{equation}
is unitary  on the boundary for a.a. $z\in\bT$. Consequently
$U(e^{it})\in L_\iy({\mathbb T})$.

Since $\chi_j^+(z)$, $j=1,2$, are  spectral factors, it is assumed
that their entries are from $H_2$ and $\det\chi_j^+(z)$, $j=1,2$,
are outer analytic functions, so that entries of
$$U(z)=\frac1{\det\chi_2^+(z)} \Adj(\chi_2^+(z))\chi_1^+(z)$$
can be represented as ratios of two functions from $H_{2/r}$ and
$H_{2/r}^O$, respectively. Hence, we can use the generalization of
Smirnov's theorem  to conclude that $U(z)\in H_\iy$, i.e.
$U(e^{it})\in L^+_\iy({\mathbb T})$.

By changing the roles of $\chi_1^+$ and $\chi^+_2$ in this
discussion, we get $$(\chi^+_1(z))^{-1} \chi_2^+(z)\in H_\iy.$$ But
$(U(z))^*=(\chi^+_1(z))^{-1}\chi_2^+(z)$ for a.a. $z\in\bT$, by
virtue of (7). Thus, we have
$$
U(e^{it})\in L^+_\iy({\mathbb T})\text{ and }\ol{U(e^{it})}\in
L^+_\iy({\mathbb T}),
$$
which implies that (8) is a constant matrix-function  and (6)
follows.
 \end{proof}

\subsection{Main Lemmas}

\begin{lemma}
For any $m\!\times\!m$ matrix-function $F_m(z)$ of the form
\begin{equation}
F_m(z)=\begin{pmatrix}1&0&0&\cdots&0&0\\
          0&1&0&\cdots&0&0\\
           0&0&1&\cdots&0&0\\
           \vdots&\vdots&\vdots&\vdots&\vdots&\vdots\\
           0&0&0&\cdots&1&0\\
           \vf_{1}(z)&\vf_{2}(z)&\vf_{3}(z)&\cdots&\vf_{m-1}(z)&f^+(z)
           \end{pmatrix},
\end{equation}
$|z|=1$, where $f^{+}\in H_2^O$ and $ \vf_j\in L_2^{N_-}({\mathbb
T})$, $j=1,2,\ldots, m-1$, for some positive integer $N$,  there
exists a unitary matrix-function $U_m(z)$ of the form
\begin{equation}
U_m(z)=\begin{pmatrix}u_{11}(z)&u_{12}(z)&\cdots&u_{1m}(z)\\
                 u_{21}(z)&u_{22}(z)&\cdots&u_{2m}(z)\\
           \vdots&\vdots&\vdots&\vdots\\
           u_{m-1,1}(z)&u_{m-1,2}(z)&\cdots&u_{m-1,m}(z)\\[2mm]
           \ol{u_{m1}(z)}&\ol{u_{m2}(z)}&\cdots&\ol{u_{mm}(z)}\\
           \end{pmatrix},
\end{equation}
$|z|=1$, where $$u_{jk}\in L_{[N]}^+\,,\;\;\; j,k=1,2,\ldots,m,$$
and
\begin{equation}
\det U_m(z)=1,
\end{equation}
$|z|=1$, such that
\begin{equation}F_m(z)\, U_m(z)\in L^{+}_2({\mathbb T}).\end{equation}
\end{lemma}

For two dimensional matrices this lemma is proved in [6], and it is
generalized for any dimensional matrices in [1]. We emphasize that
the lemma can be proved without any reference to the existence
theorem of matrix spectral factorization. Furthermore, a system of
linear equations which provides the coefficients of functions
$u_{jk}$, $j,k=1,2,\ldots,m,$ whenever $N$ negative coefficients of
$\vf_j$, $j=1,2,\ldots,m-1,$ and $N$ positive coefficients of $f^+$
are given, can be written and solved explicitly (see [1], p. 22).

\begin{lemma}
For any $m\!\times\!m$ matrix-function $F_m(z)$, $|z|=1$, of the
form $(9)$, where $f^{+}\in H_2^O$ and $ \vf_j\in L_2({\mathbb T})$,
$j=1,2,\ldots, m-1$, there exists a unitary matrix-function $U_m(z)$
of the form $(10)$ satisfying $(11)$ a.e. such that
\begin{equation}
u_{jk}(z)\in L_{\iy}^+({\mathbb T}), \;\;j,k=1,2,\ldots,m,
\end{equation}
 and $(12)$  holds.
\end{lemma}

\begin{proof}
Let $F_m^{(N)}(z)$ be the $L_2$-approximation of matrix-function
$F_m(z)$ where the entries $ \vf_j(z)$, $j=1,2,\ldots, m-1$, are
approximated by there Fourier series
$$
\vf_{j}(z)\approx\sum_{n=-N}^\iy c_n(\vf_{j})z^n.
$$
Then we can use Lemma 1 which provides the existence of unitary
matrix-function $U_m^{(N)}(z)$,
\begin{equation}
 U_m^{(N)}(z)(U_m^{(N)}(z))^*=I_m\; \text{ (a.e.)},
\end{equation}
such that
\begin{equation}
U_m^{(N)}(z)=\begin{pmatrix}u^{(N)}_{11}(z)&u^{(N)}_{12}(z)&\cdots&u^{(N)}_{1m}(z)\\[1mm]
                 u^{(N)}_{21}(z)&u^{(N)}_{22}(z)&\cdots&u^{(N)}_{2m}(z)\\
           \vdots&\vdots&\vdots&\vdots\\
           u^{(N)}_{m-1,1}(z)&u^{(N)}_{m-1,2}(z)&\cdots&u^{(N)}_{m-1,m}(z)\\[2mm]
           \ol{u^{(N)}_{m1}(z)}&\ol{u^{(N)}_{m2}(z)}&\cdots&\ol{u^{(N)}_{mm}(z)}\\
           \end{pmatrix},\;\;\;u^{(N)}_{jk}\in L_{[N]}^+,
\end{equation}
\begin{equation}
\det U_m^{(N)}(z)=1\; \text{ (a.e.)},
\end{equation}
 and
\begin{equation}F_m^{(N)}(z)\, U_m^{(N)}(z)\in L^{+}_2({\mathbb T}).\end{equation}

Now, a convergent subsequence can be extracted from
$\{U_m^{(N)}\}_{N=1}^\iy$. Furthermore, if we require in addition,
say, positive definiteness  of the matrix $F_m^{(N)} U_m^{(N)}(0)$
which can be achieved by multiplying, if necessary, $U_m^{(N)}(z)$
from the right by a constant unitary matrix with determinant $1$,
then $U_m^{(N)}(z)$ itself converges at least in measure,
\begin{equation}
U_m^{(N)}(z)\rightrightarrows U_m(z),
\end{equation}
as $N\to\iy$. These facts were proved in [6] for two dimensional
case and, in the similar way, this can be done for any dimensional
matrices as soon as the explicit form of $U_m^{(N)}(z)$ is obtained.
Anyway, the convergence (18), together with boundedness of  unitary
matrix-functions, guaranties that we can pass to the limit in
(14)-(17), so that unitary matrix-function $U_m(z)$ in (18)
satisfies the desired conditions (10)-(13).
\end{proof}

\subsection{Proof of the theorem} First  perform the lower-upper triangular factorization of
(1) with positive entries on the diagonal, i.e. take
\begin{equation}
S(z)=A(z)(A(z))^* \text{ for a.a. }z\in \bT,
\end{equation}
where
$$
A(z)=\begin{pmatrix} f_{11}(z)& 0& \cdots&0\\
f_{21}(z)& f_{22}(z)& \cdots &0\\
\vdots&\vdots&\vdots&\vdots\\f_{r1}(z)& f_{r2}(z)&
\cdots&f_{rr}(z)\end{pmatrix}
$$
and $f_{jj}(z)\geq 0$ for a.a. $z\in \bT$, $j=1,2,\ldots,r$. This
factorization can be achieved pointwise by the well known theorem of
Linear Algebra.

Since $S_{jj}(z)\in L_1 (\bT$), $j=1,2,\ldots,r$, by virtue of
equations $$\sum_{k=1}^j |f_{jk}(z)|^2=S_{jj}(z)$$ (see (19)), all
the entries of $A(z)$ are square integrable, $$f_{jk}(z)\in L_2
(\bT),$$ $1\leq j\leq r$, $1\leq k\leq j$. Furthermore, since
$$
\sum_{j=1}^r\log f_{jj}(z)=\log \prod _{j=1}^r f_{jj}(z)=\log\det
A(z)=\frac 12 \log\det S(z)\in L_1 (\bT)
$$
(see (2)), we have
\begin{equation}
\log f_{jj}(z)\in L_1 (\bT),\;\;\;\;\;\;\;\; j=1,2,\ldots,r,
\end{equation}
The condition (20) provides that
$$
f_j^+(z)=\exp\left(\frac 1{2\pi}
\int\nolimits_0^{2\pi}\frac{e^{it}+z}{e^{it}-z}\log
f_{jj}(e^{it})\,dt\right), \;\;\; j=1,2,\ldots,r.
$$
are outer analytic functions satisfying $|f_j^+(z)|=f_{jj}(z)$ for
a.a. $z\in \bT$.

If we denote the ratio $f_j^+(z)/f_{jj}(z)$ by $u_j(z)$, then
$|u_j(z)|=1$ for a.a. $z\in \bT$, $j=1,2,\ldots,r$, and $U(z)=\diag
(u_1(z),u_2(z),\ldots,u_r(z))$ is unitary matrix-function. Thus, for
matrix-function $M(z)=A(z)U(z)$, we have
\begin{equation}
S(z)=M(z)(M(z))^*  \text{ for a.a. }z\in \bT,
\end{equation}
where
\begin{equation}
M(z)=\begin{pmatrix}f^+_1(z)&0&\cdots&0&0\\
        \vf_{21}(z)&f^+_2(z)&\cdots&0&0\\
        \vdots&\vdots&\vdots&\vdots&\vdots\\
        \vf_{r-1,1}(z)&\vf_{r-1,2}(z)&\cdots&f^+_{r-1}(z)&0\\
        \vf_{r1}(z)&\vf_{r2}(z)&\cdots&\vf_{r,r-1}(z)&f^+_r(z)
        \end{pmatrix},
\end{equation}
$f_j^+(z)\in H_2^O$, $\vf_{jk}\in L_2(\bT)$, $1\leq j\leq r$, $1\leq
k\leq j$. Define
\begin{equation}
M_1(z)=M(z)\text{ and } M_m(z)=M_{m-1}(z) V_m(z),\;m=2,3,\ldots,r,
\end{equation}
where $V_m(z)$, $m=2,3,\ldots,r$, are unitary matrix-functions
\begin{equation}
V_m(z)V^*_m(z)=I_m,  \text{ for a.a. }z\in \bT,
\end{equation}
constructed recurrently as follows: let $F_m(z)$ be the
$m\!\times\!m$ matrix of the form (9), where its last row coincides
with the last row of $M_{m-1}^{[m,m]}(z)$, and let $U_m(z)$ be the
corresponding unitary matrix-function which existence is proved in
Lemma 2, so that (11) and (12) hold. Define $V_m(z)$
 as the block matrix-function
\begin{equation}
V_m(z)=\begin{pmatrix}U_m(z)&{\mathbf 0}\\{\mathbf
0}&I_{r-m}\end{pmatrix},
\end{equation}
$m=2,3,\ldots,r$. It is  assumed that $V_r(z)=U_r(z)$. Obviously
$V_m(z)$ is unitary matrix-function and (see (11))
\begin{equation}
\det V_m(z)=\det U_m(z)=1  \text{ for a.a. }z\in \bT.
\end{equation}

Pay attention that  the following equation holds
\begin{equation}
\begin{pmatrix}M_{m-1}^{[m-1,m-1]}(z)&{\mathbf 0}\\{\mathbf 0}&1\end{pmatrix}\cdot
F_m(z)=M_{m-1}^{[m,m]}(z)\,,
\end{equation}
 while
\begin{equation}
M_{m-1}^{[m-1,m-1]}(z)\in L_2^+(\bT)
\end{equation}
for $m=2$ (see (23), (22)) and for $m>2$ as well, because of the
construction process (see (30) below). Indeed, by virtue of (23) and
(25), we have
\begin{equation}
M_m^{[m,m]}(z)=M_{m-1}^{[m,m]}(z)\,U_m(z),\;\;\;m=2,3,\ldots,r.
\end{equation}
Thus, it follows from (29),  (27), (28), and (12) that
\begin{equation}
M_{m}^{[m,m]}(z)=\begin{pmatrix}M_{m-1}^{[m-1,m-1]}(z)&{\mathbf
0}\\{\mathbf 0}&1\end{pmatrix}F_m(z)U_m(z)\in L_2^+(\bT).
\end{equation}
For $m=r$, we have
\begin{equation}M_r(z)\in L^{+}_2({\mathbb T}),\end{equation}
and we conclude that $M_r(z)$ is a spectral factor of $S(z)$,
\begin{equation}\chi^+(z)=M_r(z).\end{equation}
Indeed, since (21), (23), and (24) hold, we have
$$
S(z)=M_r(z)(M_r(z))^* \text{ for a.a. }z\in \bT,
$$
and it remains to show that $\det M_r(z)$ is outer, where $M_r(z)$
is assumed extended in the unit disk $D$ in this case, by virtue of
(31).

The equations in (23) imply that
\begin{equation}
M_r(z)=M(z)U_2(z)U_3(z)\cdot\ldots U_r(z) \text{ for a.a. }z\in \bT.
\end{equation}
Hence, taking into account (33), (22), and (29), we have
\begin{equation}
\det M_r(z)= f_1^+(z) f_2^+(z)\ldots f_r^+(z) \text{ for a.a. }z\in
\bT,
\end{equation}
 The both sides of (34) are
functions from $H_{2/r}$ and they coincide on the boundary almost
everywhere. Thus the equation in (34) is valid inside the unit
circle for each $z\in D$, and since each $f_j^+(z)$ is outer,
$j=1,2,\ldots,r$, their product $\det M_r(z)$ is outer as well.

The proof of the relation (32) is completed.

\vskip+0.7cm

Authors address:

A. Razmadze Mathematical Institute

1, Aleksidze str. Tbilisi, 0193

Georgia

{\em E-mails:} lasha@rmi.acnet.ge;\;edem@rmi.acnet.ge

\end{document}